\documentclass{article}

\usepackage[a4paper,
            bindingoffset=0.2in,
            left=1in,
            right=1in,
            top=1in,
            bottom=1in,
            footskip=.25in]{geometry}

\usepackage[maxbibnames=9]{biblatex}
\usepackage{graphicx} 
\usepackage{amsmath}
\usepackage{amsfonts}
\usepackage{amssymb}
\usepackage{amsthm}
\usepackage{hyperref}
\usepackage{tikz}
\usepackage{caption}
\usepackage{orcidlink}
\usepackage{authblk}
\usepackage{lineno}

\usepackage{xcolor}
\definecolor{turchese}{RGB}{35, 174, 163}

\renewcommand{\P}{\mathcal{P}}
\newcommand{\R}{\mathbb{R}}
\newcommand{\de}{\mathrm{d}}
\newcommand{\Skn}{\mathfrak{S}_k^n}

\newtheorem{theorem}{Theorem}
\newtheorem{remark}[theorem]{Remark}
\newtheorem{proposition}[theorem]{Proposition}

\title{A least squares approach to Whitney forms}

\author[a,b]{Ludovico Bruni Bruno\footnote{email \href{mailto:ludovico.brunibruno@unipd.it}{ludovico.brunibruno@unipd.it}} \orcidlink{0000-0002-5246-8049}}

\affil[a]{\footnotesize Dipartimento di Matematica \textquotedblleft Tullio Levi-Civita\textquotedblright, Università di Padova, via Trieste, 63, Padova, 35131, Italia}

\affil[b]{\footnotesize Istituto Nazionale di Alta Matematica \textquotedblleft Francesco Severi\textquotedblright, Piazzale Aldo Moro, 5, Roma, 00185, Italia}
            
\author[c]{Giacomo Elefante\footnote{email \href{mailto:giacomo.elefante@unito.it}{giacomo.elefante@unito.it}} \orcidlink{0000-0001-5576-6802}}

\affil[c]{\footnotesize Dipartimento di Matematica \textquotedblleft Giuseppe Peano\textquotedblright, Università di Torino, via Carlo Alberto, 10, Torino, 10123, Italia}
            
\date{}

\begin{document}

\maketitle

\begin{abstract}
    In this work we describe and test the construction of least squares Whitney forms based on weights. If, on the one hand, the relevance of such a family of differential forms is nowadays clear in numerical analysis, on the other hand the selection of performing sets of supports (hence of weights) for projecting onto high order Whitney forms turns often to be a rough task. As an account of this, it is worth mentioning that Runge-like phenomena have been observed but still not resolved completely. We hence move away from sharp results on unisolvence and consider a least squares approach, obtaining results that are consistent with the nodal literature and making some steps towards the resolution of the aforementioned Runge phenomenon for high order Whitney forms.
\end{abstract}


\section{Introduction}

Interpolation and approximation are crucial aspects of the numerical analysis of both physical- and engineering-oriented problems \cite{PletzerFillmore}. These techniques make it possible to translate observations and measurements into mathematical objects on which computations may be performed. They are hence a ground step of all numerical methods \cite{ErnFEM,HairerI,Schoenberg}.

Interpolating a scalar function means representing it by a simpler object, 
typically a polynomial, 
asking that they share the same values (measurements) on a set of nodes. This is known as Lagrange interpolation \cite[Chapter $ 2 $]{Davis} and applies to physical quantities that have a natural pointwise meaning, such as the temperature. The reliability of this process is essential and in this direction much work has been done to avoid undesired behaviours such as the Runge phenomenon \cite{Runge}. 

When vector fields (e.g., forces) or tensor fields (e.g., metric tensors) come at play the situation becomes more complicated.  
These quantities are measured via diffused data, such as circulations along lines and fluxes across surfaces, and one hence has to replace nodal evaluations by \emph{currents} \cite[Chapter $ 6 $]{SimonGMT}. 
This is formalised in the mathematical language of differential forms \cite{Flanders}, which one may intuitively think of as objects that can be integrated. The domain of integration must be consistent with the order of the form: for example, the electric field has a well defined line integral, whereas the magnetic field has a meaningful surface integral \cite{BottInteractions}, unveiling the intrinsic difference between such objects. 
To retrieve a concrete extension of Lagrange interpolation to differential forms, 
%
one may invoke the theory of \emph{weights} \cite{Rapetti07,RapettiBossavit}. 

Weight-based interpolation 
consists in associating a differential form in an infinite dimensional space (such as that of smooth forms $ \Lambda^k $, see \cite{AbrahamMarsden}) with a differential form in a finite dimensional space (such as that of high order Whitney forms $ \mathcal{P}_r^- \Lambda^k $, see \cite{HiptmairHOWF}) by matching their integrals 
\cite{ABR20,ABR23}. The main drawback of this approach, which is inherited from multivariate Lagrange interpolation, lies in the complexity of determining whether a specific collection of supports induces a well defined interpolator. Sharp proofs in this direction are often heavy and confined to very specific cases \cite{BruniThesis}; this opens the gap between redundance \cite{ChristiansenRapetti} and minimality \cite{ABRenu}. 

In place of exact interpolation, one may consider \emph{least squares} approximation \cite[Chapter $ 2 $]{Rivlin}. Instead of matching the degrees of freedom, one looks for the object that minimises a distance (typically the $2$-norm) from the observations. This is a classical technique in the approximation of functions \cite{CheneyKincaid,Isaacson,TrefethenNLA}. In the context of  differential forms it found applications to fluids \cite{TonnonHiptmair} and electromagnetism \cite{LSFIM}. In this work we pursue this approach, giving a general way for constructing a least squares framework for high order Whitney forms. We observe interesting behaviours:
\begin{itemize}
    \item if the total polynomial degree is fixed and hence the redundant and the minimal sets under consideration are chosen as in \cite{ChristiansenRapetti} and \cite{ABRenu}, performances are very similar between least squares approximation and interpolation. This happens both with an easy to interpolate differential forms (see Fig. \ref{fig:chooseR}, left) and with Runge-type differential forms, such as those investigated in \cite{BruniRunge};
    \item if the total polynomial degree is fixed, increasing the number of observations over a certain quantity (depending on the degree itself) does not improve the least squares approximation error, which stabilises below the interpolation error (see Fig. \ref{fig:1DRunge}, left). This is consistent with nodal observations \cite{BoydXu}; as a consequence, this is not a possible solution to the Runge phenomenon;
    \item combining the above strategies, carefully increasing the number of observations along with the polynomial degree leads also to convergence of Runge-like forms even for \emph{uniform} (in the sense of \cite{BruniRunge}) distribution of supports (see Fig. \ref{fig:1DRunge}, right), thus defeating this undesirable effect.
\end{itemize}

The paper is organised as follows. In Section \ref{sect:WF}, we recall the main features of Whitney forms and an algorithm for the computation of bases. We introduce the concept of weight and the relative interpolation and least squares approximation problems. In Section \ref{sect:computation}, we translate the problem to the more common language of vector calculus and explain how usual quadrature formulae can be applied in this context. Section \ref{sect:numerical} contains numerical experiments and verification. In Section \ref{sect:conclusions}, we gather conclusions and overview possible further developments.

\section{Approximation by Whitney forms} \label{sect:WF}

Differential forms are nowadays a widely spread language in physics and numerical methods, see e.g. \cite{AFW,Bossavit,CotterSW,HiemstraJCP,HiptmairCanonical,PletzerFillmore,WarnickProc}. A peculiar finite dimensional space that emerges in applications is that of Whitney forms \cite{Whitney}, also known as \emph{trimmed polynomial differential forms}. They play a relevant role in computational electromagnetism \cite{Bossavit} (see also the whole collection of Bossavit's \emph{Japanese papers}, such as \cite{Bossavit4}) and, more generally, in the construction of families of finite elements \cite{AFW,HiptmairCanonical}. In particular, they parametrise (in the sense of \cite{Gopalakrishnan}) the first family of N\'ed\'elec finite elements \cite{NedelecFirst}. What makes Whitney forms special is their deep relationship with the underlying geometry of the element.

\subsection{From low-order to high-order Whitney forms}

Let $ T \subset \R^n $ be an $n$-simplex spanned by vertices $ \{ \boldsymbol{v}_i \}_{i=0}^n $ and endowed with barycentric coordinates $ \{ \lambda_i \}_{i=0}^n $, i.e., the collection of affine polynomials such that $ \lambda_i (\boldsymbol{v}_j) = \delta_{i,j} $. We denote by $ F_\sigma $ be the $k$-subsimplex of $ T $ associated with the increasing mapping $ \sigma: \{ 0, \ldots, k \} \to \{0, \ldots, n \} $, that is the $k$-simplex generated by $ \{ \boldsymbol{v}_{\sigma(0)}, \ldots, \boldsymbol{v}_{\sigma(k)} \}$, see e.g. \cite[Chapter $1$]{Prasolov}. If we let $ \mathfrak{S}_k^n \doteq \left\{ \sigma: \{ 0, \ldots, k \} \to \{0, \ldots, n \} \right\} $ denote the set of all above increasing mappings, a basis for degree $ 1 $ (or lowest order) Whitney forms is given by the \textquotedblleft Lagrange basis\textquotedblright\ dual to integration over $ F_\sigma $, as $ \sigma $ ranges over $ \Skn $
\begin{equation} \label{eq:dualWF}
\omega_{\sigma}: \quad \int_{F_{\sigma'}} \omega_{\sigma} = \delta_{\sigma, \sigma'} , 
\end{equation}
possibly up to a rescaling factor $ k! $, see \cite{ChristiansenRapetti,Whitney}. An explicit representation of Whitney forms is thus
\begin{equation} \label{eq:defWF}
    \omega_\sigma = \sum_{i=0}^k (-1)^i \lambda_{\sigma(i)} \de \lambda_{\sigma(0)} \wedge \ldots \wedge \widehat{\de \lambda}_{\sigma(i)} \wedge \ldots \wedge \de \lambda_{\sigma(k)},
\end{equation}
the hat denoting that the underlying term has been removed. We set 
$$ \P_1^- \Lambda^k (T) \doteq \mathrm{span} \left\{ \omega_\sigma \mid \sigma \in \mathfrak{S}_k^n \right\} ,$$
and from \eqref{eq:dualWF} we deduce that the dimension of this space is $ \binom{n+1}{k+1} $. This implies that $ \P_1^- \Lambda^k (T) $ is a proper subspace of the space of all degree $ 1 $ polynomial differential forms $ \P_1 \Lambda^k (T) $ unless $ k = 0 $.

There are several ways to extend the definition of the above space to $ r > 1 $, see e.g. \cite[Chapter $ 3 $]{AFW}, \cite[Section $ 3 $]{ChristiansenRapetti} and \cite[Section $ 3 $]{HiptmairCanonical}. If $ \P_r \Lambda^k (T) $ denotes the space of $k$-forms whose coefficients are $n$-variate degree $r$ polynomials, we may think of the space $ \mathcal{P}_r^- \Lambda^k (T) $ as the image of the product $ \mathbb{P}_{r-1} (T) \times \P_1^- \Lambda^k (T) \to \P_r \Lambda^k (T) $, namely
\begin{equation} \label{eq:defPR-}
\mathcal{P}_r^- \Lambda^k (T) = \left\{ \eta \in \Lambda^k (T) \mid \eta = p \cdot \omega, \; p \in \mathbb{P}_{r-1} (T), \ \omega \in \P_1^- \Lambda^k (T) \right\} .
\end{equation}
Since the dimensional inequality
\begin{equation} \label{eq:dimPR-}
\dim \mathcal{P}_r^- \Lambda^k (T) = \binom{r+k-1}{k} \binom{n+r}{n-k} \leq \binom{n+r}{n} \binom{n}{k} = \dim \mathcal{P}_r \Lambda^k (T)
\end{equation}
holds, as proved in \cite[Eq. $(61)$]{ChristiansenRapetti}, the inclusion $ \P_r^- \Lambda^k (T) \subseteq \P_r \Lambda^k (T) $ is strict unless $ k = 0 $.

The homogenization formula for a degree $ r $ polynomial allows to write a basis of $ \mathbb{P}_{r-1} (T) $ in terms of monomials in barycentric coordinates and multi-indices $ \mathcal{I} (n+1, r-1) $ of length $ n +1 $ and weight $ r-1 $, that means $ \boldsymbol{\lambda}^{\boldsymbol{\alpha}} \doteq \prod_{i=0}^n \lambda_i^{\alpha_i} $, with $ | \boldsymbol{\alpha} | \doteq \sum_{i=0}^n |\alpha_i | = r-1 $. Merging this with \eqref{eq:defWF} and \eqref{eq:defPR-}, it follows that $ \mathcal{P}_r^- \Lambda^k (T) $ is generated by forms $ \boldsymbol{\lambda}^{\boldsymbol{\alpha}} \omega_\sigma $, with $ \boldsymbol{\alpha} \in \mathcal{I} (n+1,r-1)$. 
Also, \eqref{eq:dimPR-} shows that this is necessarily a system of generators and not a basis for $ \mathcal{P}_{r}^- \Lambda^k (T) $ for $ k > 0 $ and $ r > 1 $. To retrieve a basis we shall invoke, for instance, \cite[Theorem $4.4$]{AFW}.

\begin{theorem} \label{thm:basis}
    For each $ \sigma \in \mathfrak{S}_k^n $, define $ \mathcal{I}_\sigma \doteq \{ \boldsymbol{\alpha} \in \mathcal{I}(n+1,r-1) \mid \alpha_i = 0 \quad \forall \,i < \sigma(0) \}$. The set
    $$ \left\{ \boldsymbol{\lambda}^{\boldsymbol{\alpha}} \omega_{\sigma} \mid \sigma \in \mathfrak{S}_{k}^n, \alpha \in \mathcal{I}_\sigma \right\} $$
    is a basis for $ \P_r^- \Lambda^k (T) $.
\end{theorem}

Extensions of Whitney forms to other geometries are also considered in literature \cite{BossShapes,GradHipt,Lohi21}.

\subsection{Weights} \label{sect:weights}

We saw in \eqref{eq:dualWF} that integrals over $k$-subsimplices of $ T $ give natural degrees of freedom for $ \P_1^- \Lambda^k (T) $. Since they are too few when $ r > 1 $, i.e. in the case of $ \P_r^- \Lambda^k (T) $, enriched collection of supports $ \mathcal{S} \doteq \{ s_1, \ldots, s_N \} \subset T $ are needed. These objects are called \emph{small simplices}\footnote{Historically, only the set $ s_{(\boldsymbol{\alpha}, \sigma)} $, here defined in the sequel, was named so \cite{RapettiBossavit}.} 
and integration of a differential $ k$-form $\omega $ over elements of the set $ \mathcal{S} $ is called a \emph{weight} \cite{Rapetti07}:
\begin{equation} \label{eq:weights}
    \int_{s_i} \omega .
\end{equation}
The collection $ \mathcal{S} $ of small simplices is \emph{unisolvent} for $ \mathcal{P}_r^- \Lambda^k (T) $ if the Vandermonde matrix
\begin{equation} \label{eq:VdMMatrix}
    V_{i,j} \doteq \int_{s_i} \omega_j,
\end{equation}
written with respect to any basis $ \{ \omega_j \}_{j=1}^N $ for $ \mathcal{P}_r^- \Lambda^k (T) $, is full rank. Having at hand an invertible Vandermonde matrix \eqref{eq:VdMMatrix} allows us to compute the Lagrange basis dual to weights associated with $ \mathcal{S} $, see \cite[Eq. $(3.8)$]{BruniThesis}, and hence to expand the interpolated as $ \Pi \omega = \sum_{s_i \in \mathcal{S}} \left( \int_{s_i} \omega \right) \omega_{s_i} $,
where $ \omega_{s_i} $ is such that $ \int_{s_j} \omega_{s_i} = \delta_{i,j} $. The coefficients of such an expansion are now measurements of $ \omega $ on $ \mathcal{S} $. We address to \cite{BonRap} for a computational-oriented discussion about the duality between bases and degrees of freedom.

Generically, any collection of simplicial supports is unisolvent. Since any $ s_i \in \mathcal{S} $ is parametrised by its $k+1$ vertices, any collection $ \mathcal{S} $ is an element of the Euclidean space $ (\mathbb{R}^{k+1})^N $. For any basis $ \{ \omega_j \}_{j=1}^N $ of $ \mathcal{P}_r^- \Lambda^k (T) $, via the Lasserre theorem \cite{Lasserre}, the determinant of the Vandermonde matrix \eqref{eq:VdMMatrix} defines a polynomial function of the vertices of the supports $ \{s_i\} $. Hence, the complement of the non-unisolvence condition $ \det V = 0 $ defines a dense open set of $ (\mathbb{R}^{k+1})^N $. As a consequence, we deduce that even a non-unisolvent family may be slightly perturbed to obtain a unisolvent one. 
%
%
%
%

Two relevant unisolvent sets are studied in literature. Following the original definition (see also \cite{RapettiBossavit}), a small $k$-simplex $ s_{(\boldsymbol{\alpha}, \sigma)} $ is the datum of a pair $ (\alpha, F_\sigma) $ and a contracting factor $ r $, so that
$$ s_{(\boldsymbol{\alpha}, \sigma)} = \frac{1}{r} \left( F_{\sigma} + [\boldsymbol{v}_0 | \ldots | \boldsymbol{v}_n] \boldsymbol{\alpha}^T \right) .$$
This construction is detailed in \cite[Section $2.2$]{BruniThesis}. 
It is immediate to observe that small $k$-simplices are completely supported in $ T $. 

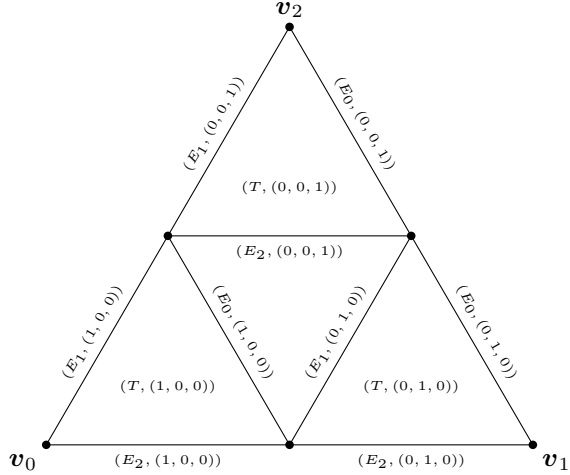
\begin{figure}[h]
	\centering
	\begin{tikzpicture}[scale = 0.8]
		\draw (0,0) -- (8,0) -- (4,6.9282) -- cycle;
		\draw (0,0) node [anchor = north east] {$ \boldsymbol{v}_0 $};
		\fill (0,0) circle (2pt);
		\draw (8,0) node [anchor = north west] {$ \boldsymbol{v}_1 $};
		\fill (8,0) circle (2pt);
		\draw (4,6.9282) node [anchor = south] {$ \boldsymbol{v}_2 $};
		\fill (4,6.9282) circle (2pt);
		\fill (4,0) circle (2pt);
		\fill (2, 3.4641) circle (2pt);
		\fill (6, 3.4641) circle (2pt);
		\draw (2, 3.4641) -- (4, 0);
		\draw (6, 3.4641) -- (4, 0);
		\draw (2, 3.4641) -- (6, 3.4641);
		\draw (2,0.7) node [anchor = south] {\tiny $ \left(T, (1,0,0) \right) $};
		\draw (6,0.7) node [anchor = south] {\tiny $ \left(T, (0,1,0) \right) $};
		\draw (4,4) node [anchor = south] {\tiny $ \left(T, (0,0,1) \right) $};
		
		\draw (2,0) node [anchor = north] {\tiny $ \left(E_2, (1,0,0) \right) $};
		\draw (3,1.73) node [anchor = south, rotate=-60] {\tiny $ \left(E_0, (1,0,0) \right) $};
		\draw (1,1.73) node [anchor = south, rotate = 60] {\tiny $ \left(E_1, (1,0,0) \right) $};
		
		\draw (6,0) node [anchor = north] {\tiny $ \left(E_2, (0,1,0) \right) $};
		\draw (7,1.73) node [anchor = south, rotate=-60] {\tiny $ \left(E_0, (0,1,0) \right) $};
		\draw (5,1.73) node [anchor = south, rotate = 60] {\tiny $ \left(E_1, (0,1,0) \right) $};
		
		\draw (4,3.4641) node [anchor = north] {\tiny $ \left(E_2, (0,0,1) \right) $};
		\draw (5,5.19) node [anchor = south, rotate=-60] {\tiny $ \left(E_0, (0,0,1) \right) $};
		\draw (3,5.19) node [anchor = south, rotate = 60] {\tiny $ \left(E_1, (0,0,1) \right) $};
	\end{tikzpicture}
	\caption{Small simplices on a triangle for $ r = 2 $.} \label{fig:smallsimplices}
\end{figure}

The following result is proved in \cite[Proposition $ 3.13 $]{ChristiansenRapetti}.

\begin{proposition} \label{prop:redundant}
    Let $ r \in \mathbb{N} $. For any $ k \geq 0 $, the set
    $$ X_r^k (T) \doteq \left\{ s_{(\boldsymbol{\alpha}, \sigma)} \mid \sigma \in \mathfrak{S}_k^n, \ \alpha \in \mathcal{I} (n+1,r-1) \right\}$$
    induces weights that determine $ \P_r^- \Lambda^k (T) $.
\end{proposition}

There is a natural bijection between the set $ X_r^k (T) $ and the system of generators $ \boldsymbol{\lambda}^{\boldsymbol{\alpha}} \omega_\sigma $, with $ \boldsymbol{\alpha} \in \mathcal{I} (n+1,r-1)$ and $ \sigma \in \mathfrak{S}_k^n $. As a consequence, $ \# X_r^k (T) \geq \dim \P_r^- \Lambda^k (T) $, and the set of Proposition \ref{prop:redundant} is minimal (hence unisolvent in a strict sense) if and only if $ k = 0, n $ or $ r = 1 $. To restore minimality one may perform the same selection as in Theorem \ref{thm:basis}. This is proved in \cite[Theorem $2.28$]{BruniThesis}.

\begin{theorem} \label{thm:minimal}
    Let $ r \in \mathbb{N} $. For any $ k \geq 0 $, the set
    $$ X_{r,min}^k (T) \doteq \left\{ s_{(\boldsymbol{\alpha}, \sigma)} \mid \sigma \in \mathfrak{S}_k^n, \ \alpha \in \mathcal{I}_\sigma (n+1,r-1) \right\}$$
    induces unisolvent weights for $ \P_r^- \Lambda^k (T) $.
\end{theorem}

The sets of Proposition \ref{prop:redundant} and Theorem \ref{thm:minimal} coincide for $ k = 0 $ and $ k = n $. In Figure \ref{fig:smallsimplices}, the difference between $ X_{2}^k (T) $ and $ X_{2,min}^k (T) $ consists only of the small edge $ \left( E_0, (1,0,0) \right) $, since $ E_0 = [\boldsymbol{v}_1, \boldsymbol{v}_2] $ is associated with the permutation $ \sigma (0,1) = (1,2) $ and $ \alpha_0 = 1 $, being $ \boldsymbol{\alpha} = (1,0,0) $. In the three dimensional framework their difference is studied in \cite{ABR20} for the case of edges $ k = 1 $ and in \cite{ABR23} for the case of faces $ k = 2 $. Despite these differences, both the above sets are relevant: the first one makes it possible to keep a rich geometrical and algebraic structure on small simplices \cite[Section $ 3 $]{ChristiansenActa}, whereas the second gives a square and invertible Vandermonde matrix \cite[Proposition $ 3.2 $]{BruniThesis}.

The possibility of geometrically choosing small simplices makes weights 
a natural extension of Lagrange interpolation to differential forms. Consistently, the choice of these supports severely affects the stability of the corresponding interpolation operator, leading to Runge-like phenomena \cite{BruniRunge}. This is explained in terms of the generalised Lebesgue constant \cite{AlonsoRapettiLeb}, although the analysis of the theoretical and numerical performances of minimal sets different from $ X_{r,min}^k (T) $ are mostly conjectural.

\subsection{Least squares approximation} \label{sect:ls}

Given the complexity of explicitly designing minimal sets of weights, in place of looking for the above described interpolated $ \Pi \omega $, one may look for a Whitney form $ \widetilde{\Pi} \omega \in \P_r^- \Lambda^k (T) $ that minimises the quantity
\begin{equation}  \label{eq:leastsquares}
    \Vert r \Vert^2 = \sum_{i=1}^M \left( \int_{s_i} \left( \widetilde{\Pi} \omega - \omega \right) \right)^2, \quad s_i \in \mathcal{S}.
\end{equation}
The resulting operator $ \widetilde{\Pi} : \Lambda^k (T) \to \mathcal{P}_r^- \Lambda^k (T) $ is a projector \cite[Chapter $ 11 $]{Powell}, i.e. $ \widetilde{\Pi} \left(\widetilde{\Pi} \omega \right)= \widetilde{\Pi} \omega $. Assume $ \mathcal{S} $ is rich and general enough, so that $ \# \mathcal{S} = M \geq N \doteq \dim \P_r^- \Lambda^k (T) $. Since $ \widetilde{\Pi} \omega \in \P_r^- \Lambda^k (T) $, it may be written as $ \widetilde{\Pi} \omega = \sum_{j=1}^{N} a_j \omega_j $, where in this case coefficients $a_j$ are not anymore measurements over $ \mathcal{S} $ but are chosen to minimise \eqref{eq:leastsquares}. 
The residual $ \Vert r \Vert^2 $ may be therefore thought as a function of variables $ \{ a_i \}_{j=1}^{N} $ which we look to minimize. Let $ W_{i,j} \doteq \int_{s_i} \omega_j $. We hence impose, for each $ j = 1, \ldots, N $,
\begin{equation} \label{eq:derivativels}
    0 = \frac{\partial \Vert r \Vert^2}{\partial a_j} = 2 \sum_{i=1}^{N} W_{i,j} \left( \sum_{\ell=1}^N W_{i,\ell} a_\ell - f_\ell \right),
\end{equation}
Up to the scaling constant, the right hand side of \eqref{eq:derivativels} is the $j$-th component of $ \left( W \boldsymbol{a} - \boldsymbol{f}\right)^T W $. As consequence, collecting the 
$N$ equations of \eqref{eq:derivativels} yields the rectangular linear system of normal equations
\begin{equation} \label{eq:derivativelssystem}
W^T W \boldsymbol{a} = W^T \boldsymbol{f} .
\end{equation}
Theorem \ref{thm:basis} comes now handy. If $ \mathcal{S} $ contains a unisolvent set (for instance considering the set described in Proposition \ref{prop:redundant}), then $ \mathrm{rank} (W) = \mathrm{rank} (W^T W) = N $ and \eqref{eq:derivativelssystem} has a unique solution
$$ \boldsymbol{a} = W^\dagger \boldsymbol{f} ,$$
where $ W^\dagger = (W^T W)^{-1} W^T $ is the Moore-Penrose pseudo-inverse of $ W $. 

\subsection{Weighted least squares}

When all the data to be fitted may be not equally relevant or reliable, least squares approximation can be appropriately modified by introducing correction coefficients. This leads to \emph{weighted least squares} \cite{Strutz} (this is an unfortunate coincidence of nomenclature: this concept of weight has nothing to do with the degrees of freedom for Whitney forms). This transforms \eqref{eq:leastsquares} into
\begin{equation}  \label{eq:weightedleastsquares}
    \Vert r \Vert^2 = \sum_{i=1}^M b_i \left( \int_{s_i} \left( \widetilde{\Pi} \omega - \omega \right) \right)^2, \quad s_i \in \mathcal{S} .
\end{equation}
Defining $ B $ as the matrix with constant diagonal entries $ B_{i,i} = b_i $, so $ B_{i,i}^{\frac{1}{2}} = \sqrt{b_i}$, one looks for
$$ \underset{\boldsymbol{a} \in \R^n}{\mathrm{arg \ min}} \left\Vert B^{\frac{1}{2}} \left( \boldsymbol{f} - W \boldsymbol{a} \right) \right\Vert .$$
The analysis may be carried exactly as that of Section \ref{sect:ls}, and normal equations \eqref{eq:derivativelssystem} read $ W^T B W \boldsymbol{a} = W^T B \boldsymbol{f} $.
If we select $ b_i = |s_i|^{-2} $, the square of the reciprocal of the measure $ |s_i| $ of the support of integration $ s_i $, we normalise the (physical) weights introduced in Section \ref{sect:weights}. This reliefs measurements from being influenced by the size of the underlying domain.

    When dealing with interpolation of Whitney forms, the norm
    \begin{equation} \label{eq:zeronorm}
        \Vert \omega \Vert_0 \doteq \sup_{c \in \mathcal{C}^k (T)} \frac{1}{|c|} \left\vert \int_c \omega \right\vert 
    \end{equation}
    is taken into account, 
    for its role in the generalised Lebesgue constant, see \cite{AlonsoRapettiLeb}. This norm is designed to provide local information, as the supremum is sought over the collection of $ k$-chains supported in $ T $, denoted by $ \mathcal{C}^k (T) $. 
    The space $ \mathcal{C}^k (T) $ is sometimes replaced by the set $\mathcal{S}^k (T)$ of simplices supported in $ T $, as in \cite{ABR20}. 
    Considering weighted least squares and observing that
    $$ \sqrt{\sum_{s_i \in \mathcal{S}} \frac{1}{|s_i|^2} \left( \int_{s_i} \omega - \widetilde{\Pi} \omega \right)^2} \geq \sup_{s_i \in \mathcal{S}} \frac{1}{|s_i|} \left\vert \int_{s_i} \omega - \widetilde{\Pi} \omega \right\vert ,$$
    we immediately deduce that a small weighted least squares approximation error over a rich collection of supports suggests a low error in the norm $ \Vert \cdot \Vert_0 $.

\subsection{Solving the system}

In practical situations, the vector $ \boldsymbol{f} $ is a collection of measurements and hence carries some error \cite{TaylorError}. The propagation of this error depends on the numerical methods used to solve the linear system, see \cite[Chapter $ 5.3 $]{Golub}. Via the SVD decomposition, one sees that the condition number $ \kappa_2 (W) \doteq \frac{\sigma_{\max} (W)}{\sigma_{\min} (W)} $ of $ W $ is a relevant quantity for the stability of the problem.

\begin{figure}[h]
    \centering
    \includegraphics[width=7.5cm]{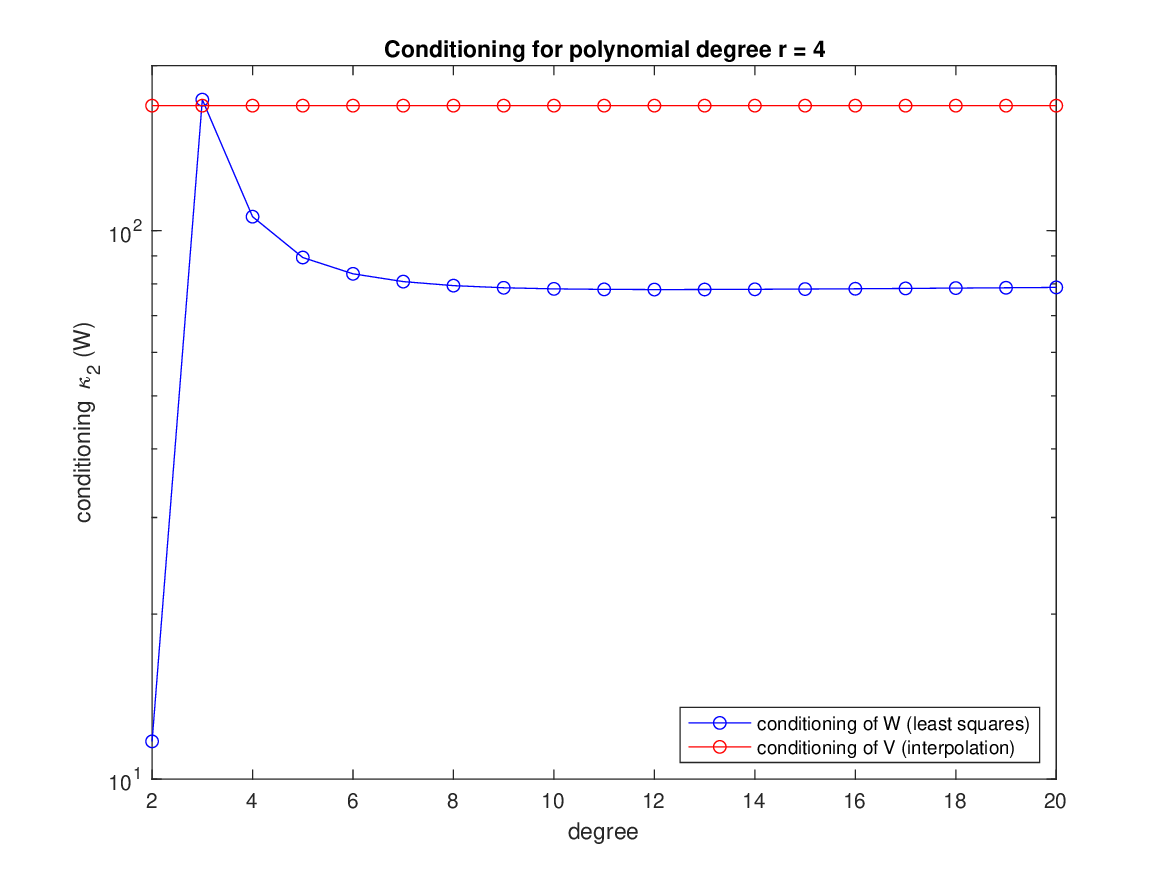}
    \includegraphics[width=7.5cm]{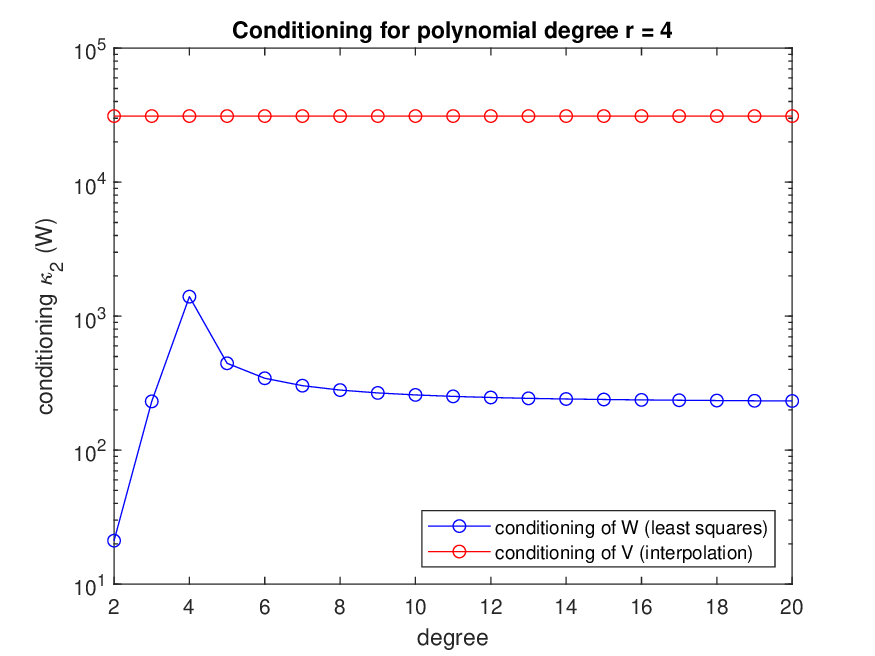}
    \caption{Conditioning of the problem for segments in $ \R^2 $ (left, $ r = 4 $) and faces in $ \R^3 $ (right, $ r = 4)$.}\label{fig:cond}
\end{figure}

Fig. \ref{fig:cond} shows an interesting behaviour of $ \kappa_2 (W) $ as a function of the supports $ X_r^k (T) $. In such a figure, the approximating polynomial degree is fixed to $ 4 $ and the set of supports $ \mathcal{S} = X_r^k (T) $ ranges in $ r = 2, \ldots, 20 $. The red line represents the conditioning of the interpolatory Vandermonde matrix $ V $. These data are consistent with the literature \cite[Tables $3.1$ and $ 3.2$]{BruniThesis}. The basis for $ \mathcal{P}_4^- \Lambda^k (T) $ is that described in Theorem \ref{thm:basis}. 
The peak of the conditioning is obtained when the least squares approximation problem is close to the interpolation problem. As $ \mathcal{S} $ is enriched, the conditioning drops and the solution remains stable. An interesting discussion about stability and accuracy of least squares problem is given in \cite{CohenLS}. 

The low conditioning 
shown in Fig. \ref{fig:cond} and the relatively small dimension of the system prescribed by Eq. \eqref{eq:dimPR-} suggest that we shall not observe a sensible difference in dependence of the method adopted for the resolution of the linear system. For ease, we hence opt for a QR factorisation. A comparison of methods adopted for solving least squares problems is given in \cite{Marchetti22}.


\section{Computation of high order weights} \label{sect:computation}

The construction of the matrix $ W $ requires the computation of weights for Whitney forms. There is a nice amount of literature about this problem. For the lowest degree Whitney forms, i.e., for elements of the space $ \P_1^- \Lambda^k (T) $, no quadrature rules are involved \cite{RapettiCompDegone}. For their high order counterparts, closed formulae may be obtained by combining localisation techniques \cite{Rapetti07} and a wise use of \textquotedblleft local\textquotedblright\ and \textquotedblleft global\textquotedblright\ barycentric coordinates \cite{ChristiansenRapetti}. Alternatively, quadrature formulae may be exploited. We expand this latter idea.


\subsection{Edge least squares} \label{sect:edgeleastsquares}

In the edge framework, the matrix $ W $ and the vector $ \boldsymbol{f} $ contain line integrals of $ 1 $-forms. 
While their meaning is clear, their computation requires a coordinate representation. 
Coordinates should be taken at convenience in order to easily handle the corresponding quadrature formulae.

Let $ T $ be an $ n $-simplex and let $ \omega \in \Lambda^1 (T) $. We have a coordinate expansion
$$ \omega = \sum_{i=1}^{n} f_i (x_1, \ldots, x_n) \de x_i .$$
Exploiting the duality $ \de x_i \left( \frac{\partial}{\partial x_j} \right) = \delta_{i,j} $ between $1 $-forms and vector fields, we may associate with $ \omega $ the vector field $ F_{\omega} = \left(f_1, \ldots, f_n \right) $. In coordinates we then have
$$ \int_{e_i} \omega = \int_{e_i} \omega \bigl|_{e_i} = \int_{e_i} j^* \omega =\int_{e_i} F_\omega \cdot \boldsymbol{t}_{e_i} \de e_i ,$$
where $ \boldsymbol{t}_{e_i} $ is the unit vector tangent to $ e_i = [\boldsymbol{p}_i, \boldsymbol{q}_i]$, i.e., $ \boldsymbol{t}_{e_i} = \frac{\boldsymbol{q}_i- \boldsymbol{p}_i}{\Vert \boldsymbol{q}_i- \boldsymbol{p}_i \Vert} $, and $ j^* : \Lambda^1 (T) \to \Lambda^1 (e_i) $ is the pullback induced by the inclusion map $ j : e_i \hookrightarrow T $. Since $ j^*\omega $ 
is a univariate function (in an appropriate chart), we may approximate the above integral with the aim of a univariate quadrature rule:
\begin{equation} \label{eq:univariatequadrature}
\int_{e_i} F_\omega \cdot \boldsymbol{t}_{e_i} \de e_i \approx \sum_{j=1}^m c_j (F_{\omega} \cdot \boldsymbol{t}_{e_i}) (\boldsymbol{x}_j) ,
\end{equation}
being $ \boldsymbol{c} = (c_1, \ldots, c_n) $ the vector of weights of the quadrature rule. The evaluation nodes $ \boldsymbol{x}_j $, for $ j = 1, \ldots, m $, belong to $ e_i $, i.e. to the segment $ \boldsymbol{p}_i + \gamma (\boldsymbol{q}_i - \boldsymbol{p}_i) $, with $ \gamma \in [0,1] $. Writing this relation componentwise, one explicitly retrieves the evaluation nodes from the one dimensional case, independently of the dimension of the ambient space. If $ m $ is chosen accordingly to the degree $ r $, the quadrature formula \eqref{eq:univariatequadrature} is exact for degree $ r $ Whitney $1$-forms, hence entries of the whole matrix $ W $ carry no approximation error.

Equation \eqref{eq:univariatequadrature} provides a way to approximate circulations. Edge least squares for Whitney forms, for polynomial degree $ r = 2 $, are invoked in literature in the context of incompressible flows \cite{TonnonHiptmair}.

\subsection{Face least squares} \label{sect:faceleastsquares}

We treat this section in the evocative framework of $2$-simplices supported in a tetrahedron $ T $. Remark \ref{rmk:codim1} extends the idea to any simplex of codimension $ 1 $ in an $ n $-dimensional framework. In a tetrahedron, the matrix $ W $ and the vector $ \boldsymbol{f} $ contain surface integrals of $2$-forms $ \omega $. In coordinates, we have
$$ \omega = f_z (x,y,z) \de x \wedge \de y + f_y (x,y,z) \de x \wedge \de z + f_x (x,y,z) \de y \wedge \de z , $$
and we want to compute
$$ \int_{s_i} \omega = \int_{s_i} \omega \bigl|_{s_i} = \int_{s_i} j^* \omega ,$$
being $ j : s_i \hookrightarrow T $ the inclusion map and $j^* : \Lambda^2 (T) \to \Lambda^2 (s_i) $ its pullback. 
We may identify the $2$-form $ \omega $ with a vector field by means of the Hodge star operator $ \star $. Letting $ F_\omega = (f_x, f_y, f_z) $, we have
$$ j^* \omega = F_\omega \cdot \boldsymbol{n}_{s_i} \de s_i ,$$
where $ \boldsymbol{n}_{s_i} $ is the outward unit normal to $ s_i $, which is well defined since $ s_i $ has codimension $ 1 $. If $ \{\boldsymbol{v}_0, \boldsymbol{v}_1, \boldsymbol{v}_2 \} $ are the vertices of $ s_i $, such a vector is readily computed as $ \boldsymbol{n}_{s_i} = \frac{(\boldsymbol{v}_1-\boldsymbol{v}_0) \times (\boldsymbol{v}_2-\boldsymbol{v}_0)}{\Vert (\boldsymbol{v}_1-\boldsymbol{v}_0) \times (\boldsymbol{v}_2-\boldsymbol{v}_0) \Vert} $. Hence, we may exploit a bivariate quadrature rule of weights $ \boldsymbol{c} = (c_1, \ldots, c_n) $:
\begin{equation} \label{eq:bivariatequadrature}
\int_{s_i} \omega = \int_{s_i} F_\omega \cdot \boldsymbol{n}_{s_i} \de s_i \approx \sum_{j=1}^m c_j (F_{\omega} \cdot \boldsymbol{n}_{s_i}) (\boldsymbol{x}_j) ,
\end{equation}
since $ j^* \omega $ is a bivariate function in an appropriate chart. 
The nodes $ \boldsymbol{x}_j $, for $ j = 1, \ldots, m $, have the same barycentric coordinates on $s_i$ that the quadrature nodes have on the reference simplex. 
Again, \eqref{eq:bivariatequadrature} is an equality when $ \omega $ is a Whitney form and $ m $ is chosen appropriately with respect to their degree $ r $, so $ W $ carries no approximation error. 

\begin{remark} \label{rmk:codim1}
    The above reasoning extends to the multi-dimensional framework as follows. First, note that an $(n-1)$-form reads in coordinates as
    $$ \omega = \sum_{i=1}^n f_i (x_1, \ldots, x_n) \de x_1 \wedge \ldots \wedge \widehat{\de x_i} \wedge \ldots \wedge \de x_n .$$
    Defining the field $ F_\omega = (f_1, \ldots, f_n) $ and computing the pullback $ j^*: \Lambda^{n-1}(T) \to \Lambda^{n-1} (s_i) $, one finds
    $$ \int_{s_i} \omega = \int_{s_i} F_\omega \cdot \boldsymbol{n}_{s_i} \de s_i ,$$
    where $ s_i $ is the $ (n-1) $-simplex spanned by vertices $ \{ \boldsymbol{v}_0, \ldots, \boldsymbol{v}_{n-1}\} $. To compute $ \boldsymbol{n}_{s_i} $ we may define $ \boldsymbol{u}_{j} \doteq \boldsymbol{v}_j - \boldsymbol{v}_0 $, and ask that $ \boldsymbol{n}_{s_i} \perp \boldsymbol{u}_j $ for $ j = 1, \ldots, n-1 $, i.e., that $ \boldsymbol{u}_j^T \cdot \boldsymbol{n}_{s_i} = 0 $. Defining the rectangular matrix $ U = \left[ \boldsymbol{u}_1^T | \ldots | \boldsymbol{u}_{n-1}^T \right] $, all normal vectors are solutions of the linear system $ U \boldsymbol{n}_{s_i} = 0 $. The desired one is obtained by normalising and imposing that $ \det \left[U^T \ | \ \boldsymbol{n}_{s_i} \right] > 0 $. 
\end{remark}

%

\section{Numerical experiments} \label{sect:numerical}
%
For computational purposes, we express \eqref{eq:univariatequadrature} and \eqref{eq:bivariatequadrature} 
in terms of the canonical basis $ \{ \boldsymbol{e}_1, \ldots, \boldsymbol{e}_n \} $.

 For $ k = 1 $, i.e., the edge case, the identification of $ 1 $-forms and vector fields given in Section \ref{sect:edgeleastsquares} is enough to obtain a Cartesian representation of Whitney forms. 
We thus get
\begin{equation}
    \de \lambda_i = \boldsymbol{e}_i \text{ if } i > 0, \quad \de \lambda_0 = - \sum_{i=1}^n \boldsymbol{e}_i .
\end{equation}

For $ k = n-1 $, to which the case of faces in a tetrahedron belongs, following the identification between $(n-1)$-forms and vector fields of Section \ref{sect:faceleastsquares} one gets
\begin{equation} \label{eq:identifn-1}
    \de \lambda_1 \wedge \ldots \wedge \widehat{\de \lambda_i} \wedge \ldots \wedge \de \lambda_n = \boldsymbol{e}_i .
\end{equation}
If $ \de \lambda_0 $ appears in the expansion, such an expression shall be replaced by $ \de \lambda_0 = - \sum_{i=1}^n \de \lambda_i $ and 
rule \eqref{eq:identifn-1} now applies, paying attention to signs that change due to alternation.

Table \ref{tab:proxies} summarises cases $ k = 1 $ in a triangle $ T \subset \R^2 $ and $ k = 2 $ in a tetrahedron $ T \subset \R^3 $. Matlab demos and codes are available at
\url{https://github.com/gelefant/LSwhitney}.

 \begin{table}[!h]
     \centering
     \begin{tabular}{c|c|| c|c}
     \multicolumn{2}{c||}{$ k = 1, \ \R^2 $} & \multicolumn{2}{c}{$ k = 2, \ \R^3$} \\
     \hline
     form & proxy & form & proxy \\
     \hline \hline
         $\de \lambda_0 $ & $ - \boldsymbol{e}_1 - \boldsymbol{e}_2 $ & $\de \lambda_0 \wedge \de \lambda_1$ & $ \boldsymbol{e}_2 + \boldsymbol{e}_3 $ \\
         $\de \lambda_1$ & $ \boldsymbol{e}_1 $ & $\de \lambda_0 \wedge \de \lambda_2$ & $ \boldsymbol{e}_1 - \boldsymbol{e}_3 $  \\
         $\de \lambda_2$ & $ \boldsymbol{e}_2 $ & $\de \lambda_0 \wedge \de \lambda_3$ & $ - \boldsymbol{e}_1 - \boldsymbol{e}_2 $  \\
         & & $\de \lambda_1 \wedge \de \lambda_2$ & $ \boldsymbol{e}_3 $\\
         & & $\de \lambda_1 \wedge \de \lambda_3$ & $ \boldsymbol{e}_2 $ \\
         & & $\de \lambda_2 \wedge \de \lambda_3$ & $ \boldsymbol{e}_1 $ \\
         \hline 
     \end{tabular}
     \caption{Explicit proxies for barycentric coordinates $1$-forms in $ \R^2 $ and $2$-forms in $ \R^3 $.}
     \label{tab:proxies}
 \end{table}

\begin{remark}
    To choose the correct degree for the quadrature formula, we observe that when the components of the field $ F_\omega $ are polynomials of degree at most $ r $, the (multivariate) function $ F_\omega \cdot \boldsymbol{n}_{s_i} $ is a (multivariate) polynomial of degree at most $ r $.
\end{remark}

\subsection{Choice of supports of approximation} \label{sect:supportselection}
%

\begin{figure}[h]
    \centering
    \includegraphics[width=7.5cm]{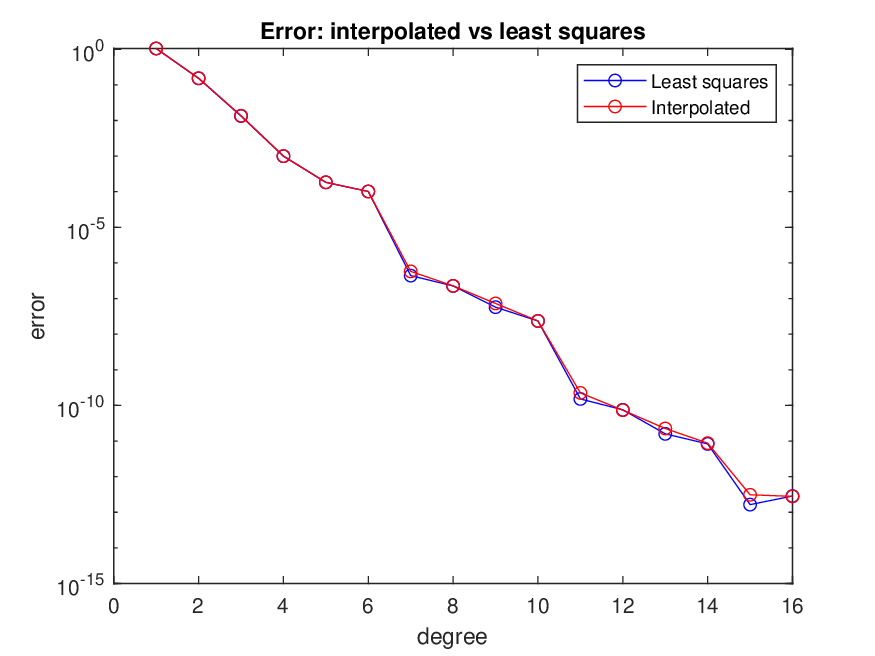} 
    \includegraphics[width=7.5cm]{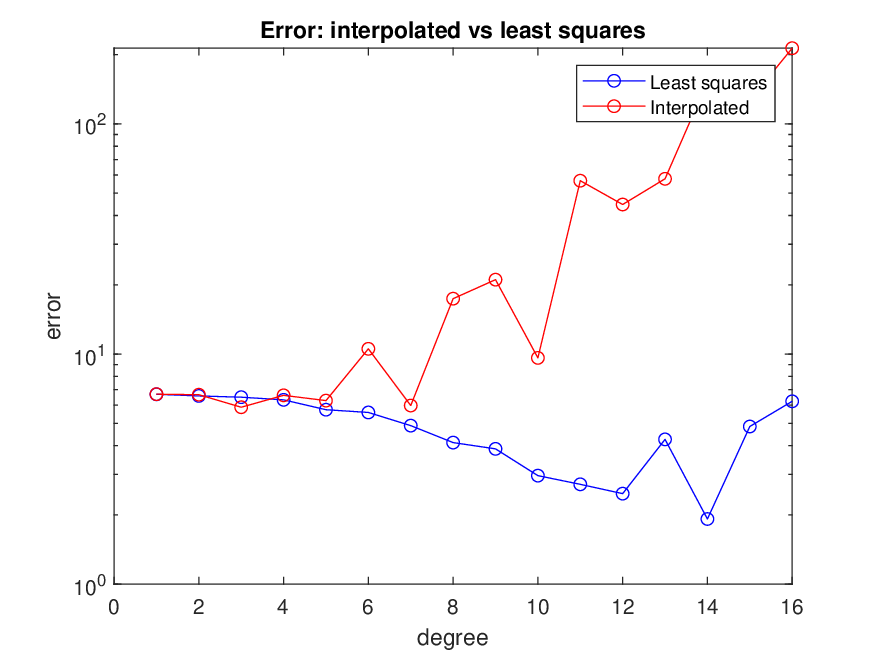}
    \caption{Left: error of the least squares approximation (blue) on $ X_r^1 (T) $ of the $1$-form $ \omega = e^{x+y} \de x + \sin{(xy)} \de y $ vs error of the interpolated (red) on minimal small simplices $ X_{r,min}^1 (T) $. Right: error of the least squares approximation (blue) of the Runge $1$-form for the set $ X^1_{r+10} (T) $ vs error of the interpolated (red) on minimal small simplices $ X_{r,min}^1 (T) $.}
    \label{fig:chooseR}
\end{figure}

When the form to reconstruct is \textquotedblleft easy\textquotedblright\ to capture (in the sense of \cite{BruniRunge}), interpolation in $ \mathcal{P}_r^- \Lambda^k (T) $ on the set $ X_{r,min}^k (T) $ and least squares approximation in $ \mathcal{P}_r^- \Lambda^k (T) $ on the set $ X_r^k (T) $ offer close results, as shown in Fig. \ref{fig:chooseR}, left. 
To improve approximation features of the least square approach we shall enrich the set of supports 
$ \mathcal{S} $ and the following question arises somewhat naturally: is the set $ X_r^k (T) $ of Proposition \ref{prop:redundant} feasible, if we replace $ r $ by a large $ R $? 
To obtain an affirmative answer,  the choice of such parameter $ R $ requires some attention. We present two possibilities:
\begin{itemize}
    \item $ R $ is significantly larger than the largest degree considered $ r_{\max} $, for instance $ R = j \cdot r_{\max} $ with $ j = 3 $: this allows to generate the set $ \mathcal{S} $ once and for all. In presence of a loop on the degree $ r $, for the lowest degrees it may give an unnecessarily large computational cost.
    \item $ R $ is a function of $ r $, for instance $ R = r + j $, with $ j = 10 $: since the set $ X_r^k (T) $ depends on $ r $, at each iteration of the loop the collection $ \mathcal{S} $ must be recomputed.
\end{itemize}

Fig. \ref{fig:chooseR} shows the sensitivity of the approximation process to this choice. If the form to be reconstructed is \textquotedblleft good\textquotedblright\, such as that defined in \eqref{eq:omega12d}, even for $ R = r $ approximation and interpolation both produce satisfactory and comparable results, see Fig. \ref{fig:chooseR}, left. In contrast, if we consider the Runge $ 1 $-form \eqref{eq:Runge2D}, the value of $ R $ must be carefully selected. In Fig. \ref{fig:chooseR}, right, we fix $ R = r + 10 $. When $ r $ is low enough (so that when we project onto $ \mathcal{P}_r^- \Lambda^k (T) $ we have $ R $ sensibly grater than $ r $), the approximated seems to converge towards the Runge form. However, as $ r $ increases, this convergence deteriorates. Errors are computed with respect to $\Vert \cdot \Vert_0 $. This example is detailed in Section \ref{sect:num2d}.

%

\subsection{On the real line}

The case $ k = 0 $ corresponds to the classical nodal least square interpolation, see e.g. \cite[Chapter 8]{Davis} and \cite[Chapter 2]{TrefethenNLA}. We consider and extend a peculiar aspect of least squares, which is treated in literature. We consider the Runge \cite{Runge} function $ f: [-1,1] \to \R $ given by
\begin{equation} \label{eq:Rungefunct}
f(x) = \frac{1}{1+25x^2} . 
\end{equation}
Such a function is well known to be not suited for interpolation at equispaced nodes. In fact, if $ r \geq 1 $ is the polynomial degree of the interpolant $ p_r (x) $ at equispaced nodes $ x_i = -1 + \frac{2i}{r} $, the quantity $ \Vert p_r (x) - f(x) \Vert_\infty $ is well known to diverge as $ r $ increases (see, e.g., \cite[Chapter 4]{Davis}).

In \cite{BoydXu}, authors try to find a least square solution to the Runge phenomenon induced by interpolation on equispaced nodes. They fix a low polynomial degree $ r $, refine the grid $ x_i = -1 + \frac{2i}{M} $ with $ M $ sufficiently large and see if the polynomial $ p_r : [-1,1] \to \R $ that solves the least square problem
\begin{equation} \label{eq:nodalLS1D}
\sum_{i=1}^M \left\Vert p_r(x_i) - \frac{1}{1+25(x_i)^2} \right\Vert_2 
\end{equation}
is in fact close enough to the Runge function \eqref{eq:Rungefunct}. The answer is negative.

\begin{figure}[h]
    \centering
    \includegraphics[width=7.5cm]{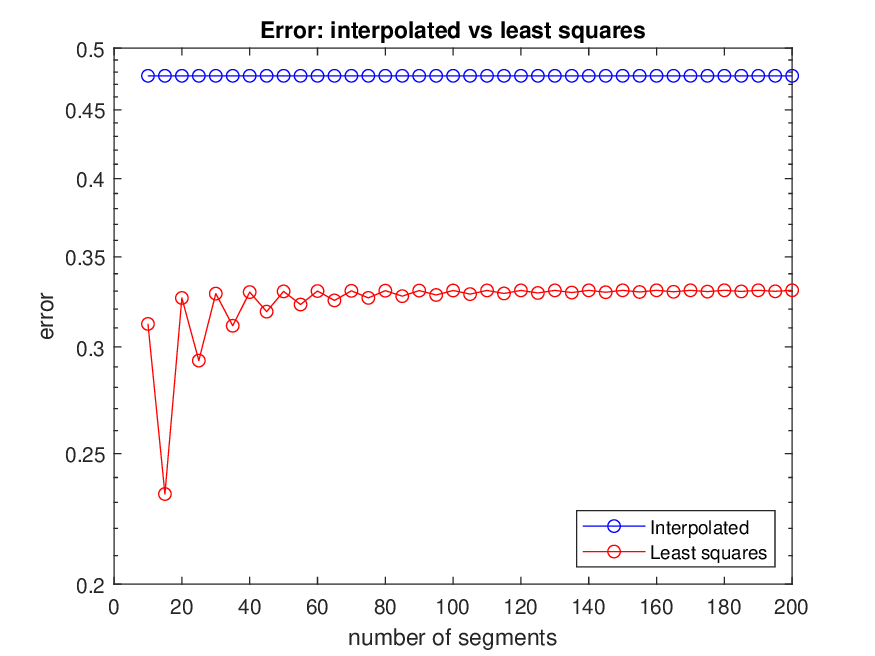}
    \includegraphics[width=7.5cm]{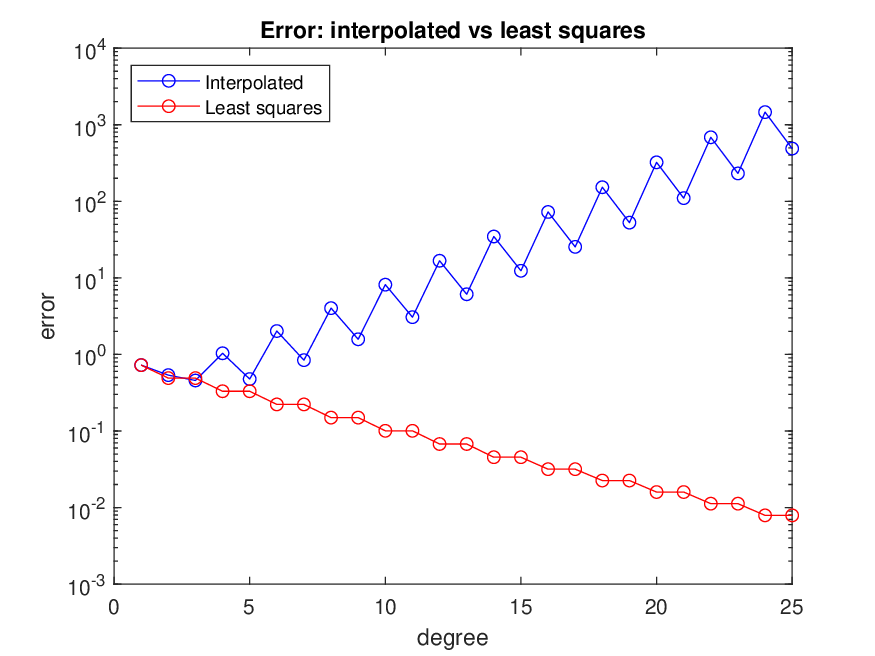}
    \caption{Left: the degree is fixed ($ r = 5$) and the number of segments is increased. Right: the number of segments is fixed (and large), the degree $ r $ is increased. Uniform segments are considered.}\label{fig:1DRunge}
\end{figure}

The case $ k = 1 $ is that of segmental interpolation (or \emph{histopolation} \cite{Robidoux}). 
We consider again the function \eqref{eq:Rungefunct} and, in place of considering nodal evaluations on equispaced points $ x_i = -1 + \frac{2 i }{r} $, we consider as degrees of freedom its \emph{uniform weights},
namely integrals of \eqref{eq:Rungefunct} on $ s_i = \left[ -1+\frac{2i}{r}, -1+\frac{2(i+1)}{r} \right] $, for $ i = 0, \ldots, r-1 $. That the corresponding interpolating polynomial do not converge towards $ f $ is known \cite{BE23}. We thus set up \eqref{eq:nodalLS1D} to this framework, considering a large $ M $, uniform segments $ s_i = \left[ -1+\frac{2i}{M}, -1+\frac{2(i+1)}{M} \right] $, and we look for a polynomial $ p_r : [-1, 1] \to \R $ of degree $ r $ such that
\begin{equation} \label{eq:segmentalLS1D}
    \sum_{i=1}^M \left\Vert \frac{1}{|s_i|} \int_{s_i} p_r(x) \de x - \frac{1}{|s_i|} \int_{s_i} \frac{1}{1+25(x_i)^2} \de x \right\Vert_2 
\end{equation}
is the least possible. Since any family of $ r + 1 $ segments in the set $ \{ s_i \}_{i=1}^M $ is unisolvent for $ \mathbb{P}_r $ (see \cite[Theorem 4.1]{Robidoux} or \cite[Proposition 3.1]{BE23}), the solution to the least squares problem \eqref{eq:segmentalLS1D} is unique. Consistently with the nodal case, the left hand panel of Figure \ref{fig:1DRunge} shows that the increment in the number of segments, when the degree is fixed (in the picture, $ r = 5 $), does not push toward convergence. To defeat the Runge phenomenon and hence obtain convergence, it is necessary to combine an increment of the degree with a large number of supports, as shown in the right hand panel of Figure \ref{fig:1DRunge}.

An explicit implementation strategy for this framework can be derived from \cite{BruniDRNA}.

\subsection{Segments in the triangle} \label{sect:num2d}

\begin{figure}[h]
    \centering
    \includegraphics[width=7.5cm]{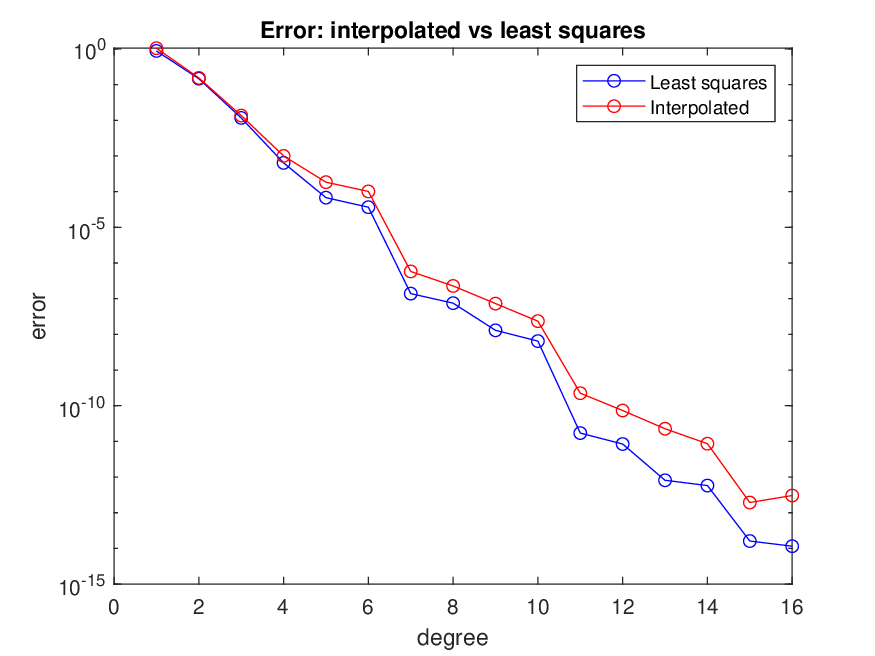}
    \includegraphics[width=7.5cm]{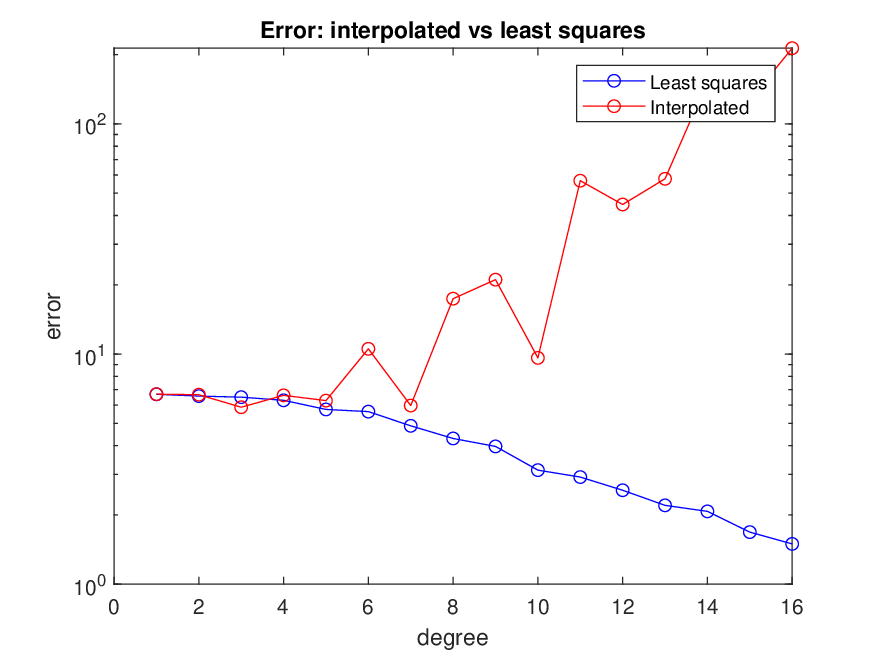}
    \caption{Zero-norm of the error with respect to the total polynomial degree. Left: results for the $1$-form $ \omega_1 $. Right: results for the Runge $1$-form.}\label{fig:2dinterp}
\end{figure}

We perform two numerical tests to test the effectiveness of the least square approximation of Whitney $1$-forms. First, we consider the form
\begin{equation} \label{eq:omega12d}
     \omega_1 \doteq e^{x+y} \de x + \sin(xy) \de y .
\end{equation}
We compute its interpolated with respect to weights of the set $ X_{r,min}^1 (T) $ of Theorem \ref{thm:minimal} and its approximated with respect to least squares over the set $ X_r^1 (T) $ of Proposition \ref{prop:redundant}. We then take a rich collection of 7056 segments supported in $ T $ and here estimate the error with respect to \eqref{eq:zeronorm}. Results of the two techniques are comparable, see Fig. \ref{fig:2dinterp}, left.

We then perform the same test on a Runge-like $1$-form. The construction of these forms is based on generalisations of the univariate Runge function \cite{Runge} to the multivariate framework \cite{RothThesis}. 
The differential of any function suffering from the Runge phenomenon is in fact a $ 1 $-form that presents the same behaviour with respect to the norm \eqref{eq:zeronorm}, see \cite[Proposition $3.4$]{BruniRunge}. A neat strategy in the choice of supports to counter this effect has not been found yet, and only mild solutions have been obtained with interpolation strategies (see \cite[Section $ 4 $]{BruniRunge} and \cite[Chapter $ 4 $]{BruniThesis}).

We follow the above strategy and consider the function
$$ f(x,y) = \frac{1}{1 + \alpha \left( x - \frac{1}{3} \right)^2} \frac{1}{1 + \alpha \left( y - \frac{1}{3} \right)^2} ,$$
where $ \alpha $ is sufficiently large to observe the phenomenon (following \cite{BruniRunge} and \cite{RothThesis}, in the experiments we take $ \alpha = 100 $) and the shift of factor $ 1/3 $ centers the function to the standard simplex. Then
\begin{equation} \label{eq:Runge2D}
\de f = \frac{2 \alpha \left(x-\frac{1}{3}\right)}{\left( 1+ \alpha \left(x-\frac{1}{3}\right)^2 \right)^2 \left( 1 + \alpha \left( y - \frac{1}{3} \right) \right)} \de x + \frac{2 \alpha \left(y-\frac{1}{3}\right)}{\left( 1+ \alpha \left(x-\frac{1}{3}\right) \right) \left( 1 + \alpha \left( y - \frac{1}{3} \right)^2 \right)^2} \de y .
\end{equation}
The complexity of finding a collection of weights that gives converge of the interpolated is discussed in \cite[Section $ 4 $]{BruniRunge} and \cite[Section $4.4$]{BruniThesis}. To be more definite, in the numerical tests available in literature, what is achieved is a non-increasing trend of the error as the degree increases. Fig. \ref{fig:2dinterp}, right, consistently with the behaviour observed on the real line, shows that a combination of least squares and redundant simplices of Proposition \ref{prop:redundant} provide a possible solution to mitigate the Runge effect. In our experiments, for any degree $ r =1,\dots,r_{\max}$, as supports for the least squares approximation we take $ X_{3r_{\max}}^1 (T) $.

\subsection{Faces in the tetrahedron}

To implement the quadrature formula \eqref{eq:bivariatequadrature} we start from a two-dimensional quadrature formula on a two-dimensional simplex. 
We embed this triangle in $ \R^3 $ by adding a third coordinate and setting it to zero. Call this $ \widehat{T} $. We shall map $ \widehat{T} $ to  any $2$-face $ T_i = \{ \boldsymbol{y}_0, \boldsymbol{y}_1, \boldsymbol{y}_2 \} \subset T \subset \R^3 $. Put
$$ \boldsymbol{v}_j = \boldsymbol{x}_j - \boldsymbol{x}_0 \text{ for } j = 1,2 , \quad \widehat{\boldsymbol{n}} = \boldsymbol{v}_1 \times \boldsymbol{v}_2 $$
and 
$$ \boldsymbol{w}_j = \boldsymbol{y}_j - \boldsymbol{y}_0 \text{ for } j = 1,2 , \quad \boldsymbol{n}_i = \boldsymbol{w}_1 \times \boldsymbol{w}_2 .$$
Since the above are both bases for $ \R^3 $, the affinity $ \varphi_i(\boldsymbol{x}) = A_i \boldsymbol{x} + \boldsymbol{b}_i $ that maps $ \widehat{T} $ to $ T_i $ (hence quadrature nodes) is completely determined by the set of equations
$$\begin{cases}
    A_i \boldsymbol{v}_j = \boldsymbol{w}_j \quad j = 1, 2 \\
    A_i \widehat{\boldsymbol{n}} = \boldsymbol{n}_i \\
    \boldsymbol{b}_i = \boldsymbol{y}_0 - A_i \boldsymbol{x}_0 .
\end{cases} $$

The matrix $ A_i $ and the vector $ \boldsymbol{b}_i $ must be determined for each triangle $ T_i $. To avoid the resolution of this linear system, we start with a quadrature formula on the standard triangle in $ \R^2 $. 
Under this hypothesis, the reference triangle reads $ \widehat{T} = \{ \boldsymbol{0}, \boldsymbol{e}_1, \boldsymbol{e}_2 \} $, whence $ \boldsymbol{v}_1 = \boldsymbol{e}_1 $, $ \boldsymbol{v}_2 = \boldsymbol{e}_2 $ and $ \widehat{\boldsymbol{n}} = \boldsymbol{e}_3 $. Thus
\begin{equation}
    A_i = [\boldsymbol{y}_1 - \boldsymbol{y}_0 \ | \ \boldsymbol{y}_2 - \boldsymbol{y}_0 \ | \ \boldsymbol{n}_i] \text{ and } \boldsymbol{b}_i = \boldsymbol{y}_0 .
\end{equation}

\begin{figure}[h]
    \centering
    \includegraphics[width=7.5cm]{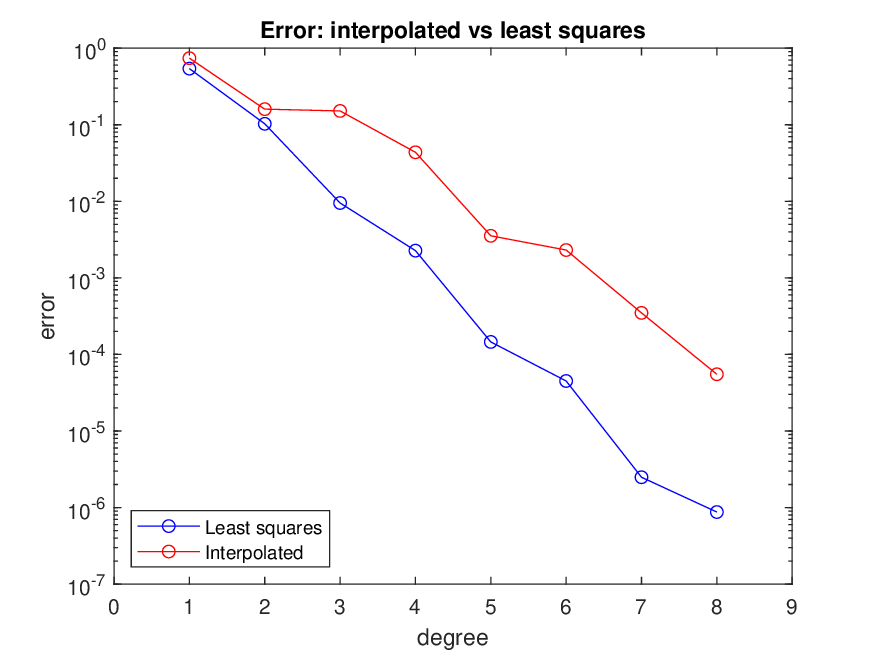}
    \includegraphics[width=7.5cm]{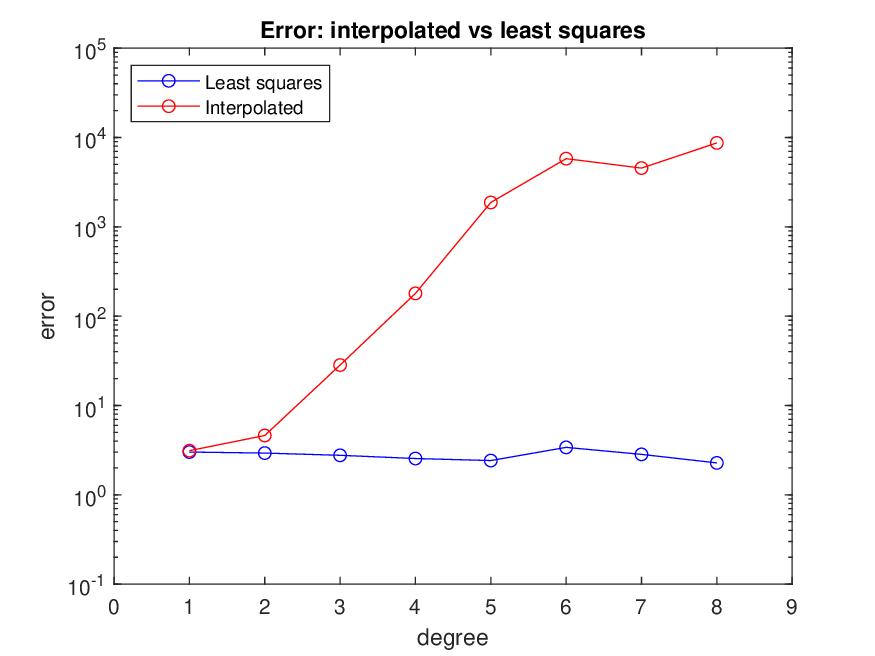}
    \caption{Zero-norm of the error with respect to the total polynomial degree. Left: results for the $2$-form $ \omega_2 $. Right: results for the Runge $2$-form.}\label{fig:3dinterp}
\end{figure}

We again perform two numerical experiments. In the first we consider the $2$-form
$$ \omega_2 \doteq e^{z+xy}\de x \wedge \de y + \sin (xy) \sqrt{z+1} \de x \wedge \de z + \cos(x+y) \de y \wedge \de z $$
and compare the interpolated with respect to the set $ X_{r,min}^2 (T) $ with the approximated obtained from the set $ X_r^2 (T) $, in terms of the $0$-norm \eqref{eq:zeronorm}. Results are close, as shown in Fig. \ref{fig:3dinterp}.

For the second case, we construct a Runge-like $2$-form. Although they are not explicit in literature, we may exploit the identification given in section \ref{sect:faceleastsquares} and take the Runge $2$-form $ \omega_2 $ whose coefficients (up to reordering) coincide with the Runge $ 1 $-form obtained by differentiating
$$ f(x,y,z) = \frac{1}{1 + \alpha \left( x - \frac{1}{4} \right)^2} \frac{1}{1 + \alpha \left( y - \frac{1}{4} \right)^2} \frac{1}{1 + \alpha \left( z - \frac{1}{4} \right)^2} .$$
The shift centres the function with the barycenter of the standard tetrahedron $ T $ and $ \alpha = 10 $, as in Section \ref{sect:num2d}. Consistently with the segmental framework, the interpolated over the set $ X_{r,min}^2 (T) $ diverges, and so does the least squares approximation obtained on $ X_r^2 (T) $. In contrast, a combination of the increment of the polynomial degree $r$, with $r=1,\dots,r_{\max}$, and a richer least square approximation, obtained over the set $ X_{r_{\max}+10}^2 (T) $ shows convergence towards the Runge $2$-form.

In all the above numerical test the collection of supports used to estimate the norm of the error contains $ 13680 $ triangles. 

\section{Conclusions and open problems} \label{sect:conclusions}

In this paper, we provided a framework for the construction of Whitney forms based on least squares weights. This bypasses the hard task of studying unisolvence of sets of $ k$-simplices. When the sets of redundant supports are entangled with the corresponding minimal sets, the behaviour of the least squares based approximant is close to the interpolant. In contrast, when these former sets are sensibly enriched, the approximated and the interpolated form diverge one from the other. If this is not particularly relevant on easy to capture objects, it becomes relevant when dealing with Runge-like forms. As a consequence, a wise combination of the total polynomial degree and the enrichment of the set supports may be used to tame the Runge phenomenon.

Some deep questions, which we aim to treat in the future, remain open. First of all, we presented interpolation and approximation as two antithetic techniques. This is not completely fair, as blended methods may be applied: for instance, we could ask a polynomial to fit some data and minimise the distance from others, as in \cite{Mazza}. This is for sure a relevant development, which requires also a non trivial understanding of the geometry of the problem.

We selected a fixed set of supports of integration in the construction of the matrix $ W $ (see discussion in Section \ref{sect:supportselection}). This may be expensive. To contain the computational cost, it is of interest to retrieve higher-degree matrices by exploiting lower-degree ones. This depends on the structure of the space $ \mathcal{P}_r^- \Lambda^k (T) $ and its decomposition. This topic is of self-standing interest, see e.g. \cite[Section $5$]{HiptmairHOWF}.

Further, in an academic problem, the error of the right hand side can be handled by considering a sufficiently strong quadrature rule and all sizes at play are controlled. In a real physical or engineering problem the number of observations might sharply grow and the vector $ \boldsymbol{f} $ may contain a non trivial noise. The reconstruction of the approximant requires a thorough study of the numerical linear algebra of the problem: stable methods such as QR or SVD factorisation bear a large computational cost; 
in contrast, faster methods might overlook the propagation of the error. For large systems, the analysis of stabilization methods and involved techniques is thus of interest.







\section*{Acknowledgements}

The authors thank Mariarosa Mazza for her precious remarks on the topic and Alvise Sommariva for his suggestions and repository of quadrature formulae on triangles available at \cite{AlviseRepository}. This research has been accomplished within Rete ITaliana di Approssimazione (RITA), the thematic group on \lq\lq Teoria dell'Ap\-pros\-si\-ma\-zio\-ne e Applicazioni\rq\rq $\,$(TAA) of the Italian Mathematical Union. The authors are members of the Gruppo Nazionale Calcolo Scientifico-Istituto Nazionale di Alta Matematica (GNCS-IN$\delta$AM), which partially supported this research. The first author is funded by IN$\delta$AM and supported by Universit\`a di Padova. The second author has been supported by Fondazione CRT, project 2022 \lq\lq Modelli matematici e algoritmi predittivi di intelligenza artificiale per la mobilit$\grave{\text{a}}$ sostenibile\rq\rq.

\printbibliography

\end{document}